\documentclass[journal]{IEEEtran}
\usepackage{amssymb}
\usepackage{graphicx}
\usepackage{titlesec}
\usepackage{amsmath}
\usepackage{diagbox}
\usepackage{xcolor}

\ifCLASSOPTIONcompsoc
\usepackage[caption=false,font=normalsize,labelfont=sf,textfont=sf]{subfig}
\else
\usepackage[caption=false,font=footnotesize]{subfig}
\fi

\def\x{\textit{\textbf{x}}}

\def\p{\textit{\textbf{p}}}

\def\n{\textit{\textbf{n}}}

\def\R{\mathbb{R}}

\def\dv{v^\prime_{_\lambda}}
\def\ddv{v^{\prime\prime}_{_{\lambda\lambda}}}
\def\us{u^\prime_{s}}
\def\uss{u^{\prime\prime}_{ss}}

\def\div{\mathop{\mathrm{div}}}

\DeclareMathOperator{\dist}{dist}

\begin{document}
\title{Accuracy Improvements for Convolutional and Differential Distance Function Approximations}

\author{Alexander~Belyaev and Pierre-Alain~Fayolle}
\maketitle

\markboth{Journal of \LaTeX\ Class Files,~Vol.~XX, No.~X, July~2024}%
{Shell \MakeLowercase{\textit{et al.}}: Bare Demo of IEEEtran.cls for IEEE Journals}

\begin{abstract}
Given a bounded domain, we deal with the problem of estimating the distance function from the internal points of the domain to the boundary of the domain. Convolutional and differential distance estimation schemes are considered and, for both the schemes, accuracy improvements are proposed and evaluated. Asymptotics of Laplace integrals and Taylor series extrapolations are used to achieve the improvements.
\end{abstract}

\IEEEpeerreviewmaketitle

%%%%%%%%%%%%%%%%%%%%%%%%%%%%%%%%%%%%%%%%%%%%%%%%%%%%%%%%%%%%%%%%%%
\section{Introduction and contribution}\label{sec:Intro}
Distance functions are widely used in analysis of signals defined on surfaces and graphs and in image analysis.
Distance function estimation and approximation remains an active research area due to applications in computational physics \cite{Park-etal_jcp24}, numerical methods for partial differential equations \cite{King-etal_tog24,MattosDaSilva-etal_tog24}, big data processing \cite{Huguet-etal_NeurIPS23}, geometric modeling \cite{Edelstein-etal_sig23,Feng-Crane_tog24,Yang-etal_NeurIPS23}, graph theory \cite{Meister_arXiv23}, and several other fields.

In this paper, we deal with the following basic distance approximation problem. Given a bounded domain $\Omega$ in $\R^{d+1}$ with smooth boundary $\partial\Omega$, our goal is to get an accurate estimation of $\dist(\x,\partial\Omega)$, the distance to $\partial\Omega$ for all points $\x\in\Omega$. Some of the distance function estimation improvements that we introduce can also be extended to more general problems including estimating geodesic distances on a surface.

Our contributions are as follows. We establish a link between convolution-based distance transforms \cite{Karam-etal_spl19} and Laplace approximations considered in \cite{Tibshirani-etal_arXiv24} in the context of non-convex variational problems. We use this connection to improve the distance function approximation accuracy by combining convolutional distance approximations. Next, we demonstrate how a convolutional distance transform scheme can be combined with the popular heat method \cite{Varadhan_cpam67,Crane-Weischedel-Wardetzky_tog13} for distance function approximation. Finally, and this is probably the main practical contribution of the paper, we propose two Taylor-based extrapolation schemes for enhancing and extending the heat distance method for flat domains.

%%%%%%%%%%%%%%%%%%%%%%%%%%%%%%%%%%%%%%%%%%%%%%%%%%%%%%%%%%%%%%%%%%
\section{Laplace's method}\label{sec:Laplace}
Let $\varphi(\p)$, $\varphi:\R^d\rightarrow\R$, be smooth and $h(\p)$, $h:\R^d\rightarrow\R$, be continuous. Let $\p^{*}$ be the global minimizer of $\varphi(\x)$, $\varphi^{*}=\varphi(\p^{*})$. Laplace's method implies that
\begin{equation}\label{eq:Laplace}
\lambda^d\!\int h(\p)e^{-\lambda\varphi(\p)}d\x\sim C\,\lambda^{d/2}e^{-\lambda\varphi(\p^{*})}h(\p^{*}),
\end{equation}
where $C$ is a positive constant and $\sim$ means that the ratio of two sides tends to $1$, as $\lambda\to\infty$. See \cite{Tibshirani-etal_arXiv24} for details.

Applying $-\frac{1}{\lambda}\log(\cdot)$ to both parts of (\ref{eq:Laplace}) yields
\begin{equation}\label{eq:underestimate}
-\frac{1}{\lambda}\log\!\left[\lambda^d\!\!\int\!\!h(\p)e^{-\lambda\varphi(\p)}d\p\!\right]\!
=\varphi^{*}-\frac{d}{2}\frac{\log\lambda}{\lambda}
+O\!\left(\frac{1}{\lambda}\right)\!,
\end{equation}
as $\lambda\to\infty$. Thus, for a sufficiently large $\lambda$, the right-hand side of (\ref{eq:underestimate}) underestimates $\varphi(\p^{*})$.

Following \cite{Tibshirani-etal_arXiv24} let us consider the following self-normalized version of (\ref{eq:Laplace})
\begin{equation}\label{eq:selfLaplace}
\min h(\p)<\frac{\int h(\p)e^{-\lambda\varphi(\p)}d\p}{\int e^{-\lambda\varphi(\p)}d\p}\sim h(\p^*).
\end{equation}
If we set $h(\p)=\varphi(\p)$, then (\ref{eq:selfLaplace}) overestimates $\varphi^{*}=\varphi(\p^{*})$. Thus
\begin{equation}\label{eq:overestimateLaplace}
\frac{\int\varphi(\p)e^{-\lambda\varphi(\p)}d\p}{\int e^{-\lambda\varphi(\p)}d\p}=\varphi^*
+\frac{K}{\lambda}+O\left(\frac{1}{\lambda^2}\right)
\end{equation}
for $\lambda\to\infty$ with constant $K>0$. We can combine estimates (\ref{eq:underestimate}) with $h(\p)\equiv1$ and (\ref{eq:overestimateLaplace}) to improve the approximation accuracy
$$
\alpha\left\{\frac{\int\varphi(\p)e^{-\lambda\varphi(\p)}\,d\p}{\int e^{-\lambda\varphi(\p)}d\x}\right\}
+\beta\left\{-\frac{1}{\lambda}\log\left[\lambda^d\!\int e^{-\lambda\varphi(\p)}d\p\right]\right\}
$$
with weights $\alpha$ and $\beta$ given by
\begin{equation}\label{eq:weights}
\alpha=\frac{d\log\lambda}{2K+d\log\lambda}\quad\mbox{and}\quad\beta=\frac{2K}{2K+d\log\lambda}.
\end{equation}
The idea to combine (\ref{eq:underestimate}) and (\ref{eq:overestimateLaplace}) seems new.
% This can be thought as a particular case of the well-known Richardson extrapolation procedure \cite{Richardson_1911}.

While constant $K$ in (\ref{eq:overestimateLaplace}) can be determined analytically,
% (see, for example, Proposition\,4.5 on p.\,125 in \cite{Fedoruk_book87} where a full asymptotics exapsion of the integral in (\ref{eq:Laplace}) is presented),
it depends on the values of $\varphi(\p)$ and its derivatives at $\p^*$ which typically are not known to us. So, in practice we have to choose $K$ heuristically.

Note that the results remain true if we consider functions $\varphi(\p)$ and $h(\p)$ defined on a $d$-dimensional closed surface (manifold). In particular, as shown in the next section, we can consider a point $\x$ inside domain $\Omega$ bounded by $\partial\Omega$ and define $\varphi(\p)$ by $\varphi(\p)=\|\x-\p\|$, where $\p$ lies on $\partial\Omega$. Then (\ref{eq:underestimate}) with $h(\p)\equiv\varphi(\p)$ and (\ref{eq:overestimateLaplace}) deliver estimates of $\dist(\x,\partial\Omega)$, the distance function to the boundary of the domain.

%%%%%%%%%%%%%%%%%%%%%%%%%%%%%%%%%%%%%%%%%%%%%%%%%%%%%%%%%%%%%%%%%%
\section{Convolutional distance function estimation}
The following two formulas are the main ingredients of fast convolutional distance transforms proposed in \cite{Karam-etal_spl19}.
\begin{equation}\label{eq:logConv}
\varphi^{*}=\min\{\varphi_1,\dots,\varphi_n\}\approx-\frac{1}{\lambda}
\log\left[\sum\limits_{k=1}^n e^{-\lambda\varphi_k}\right],
\end{equation}
\begin{equation}\label{eq:softMin}
\varphi^{*}=\min\{\varphi_1,\dots,\varphi_n\}\approx\frac{\sum_{k=1}^n\varphi_k e^{-\lambda\varphi_k}}{\sum_{k=1}^n e^{-\lambda\varphi_k}},
\end{equation}
as $\lambda\to\infty$.
One can observe that (\ref{eq:underestimate}) with $h(p)\equiv1$ and (\ref{eq:overestimateLaplace}) are continuous counterparts of (\ref{eq:logConv}) and (\ref{eq:softMin}), respectively. In particular, it is easy to see that the right-hand side of (\ref{eq:logConv}) underestimates $\varphi^{*}$,  while the right-hand side of (\ref{eq:softMin}) overestimates it.

Given a domain $\Omega$ in $\R^d$ bounded by $\partial\Omega$, let us set $\varphi(\p)=\|\x-\p\|$. Then
the continuous versions of the discrete convolutions (\ref{eq:logConv}) and (\ref{eq:softMin}) used in \cite{Karam-etal_spl19} are given by the left-hand and right-hand sides of
$$
-\frac{1}{\lambda}\log\left[\lambda^d\!\int_{\partial\Omega}e^{-\lambda\varphi(\p)}d\p\right]
\approx\varphi^{*}\approx
\frac{\int_{\partial\Omega}\varphi(\p)e^{-\lambda\varphi(\p)}d\p}{\int_{\partial\Omega}e^{-\lambda\varphi(\p)}d\p},
$$
respectively. Thus, in view of (\ref{eq:underestimate}) and (\ref{eq:overestimateLaplace}), we can improve the approximation accuracy of (\ref{eq:logConv}) and (\ref{eq:softMin}) by blending their right-hand sides with the weights defined by (\ref{eq:weights}).

As pointed out in \cite{Karam-etal_spl19}, an equivalent way to represent the right-hand side of (\ref{eq:softMin}) consists of differentiating the logarithm in the right-hand side of (\ref{eq:logConv}):
\begin{equation}\label{eq:dLogConv}
-\frac{\partial}{\partial\lambda}\log\left[\sum\limits_{k=1}^n e^{-\lambda\varphi_k}\right]
=\frac{\sum_{k=1}^n\varphi_k e^{-\lambda\varphi_k}}{\sum_{k=1}^n e^{-\lambda\varphi_k}}.
\end{equation}
This identity relating (\ref{eq:logConv}) and (\ref{eq:softMin}) looks trivial, but in the next section we will use it to get new formulas to estimate the distance function.

%%%%%%%%%%%%%%%%%%%%%%%%%%%%%%%%%%%%%%%%%%%%%%%%%%%%%%%%%%%%%%%%%%
\section{Differential distance function estimation}
A significant progress in distance function estimation was achieved in \cite{Crane-Weischedel-Wardetzky_tog13} (see also \cite{Crane-Weischedel-Wardetzky_cacm17}) where the so-called heat method was introduced. The heat method is typically used for geodesic distance estimation and based on S.\,R.\,S.\,Varadhan's classic result in differential geometry \cite{Varadhan_cpam67}. Below we consider a simplified version of the heat method adapted for planar distance function estimation.

Consider the following screened Poisson equation in a bounded domain $\Omega$
\begin{equation}\label{eq:screendPoisson}
-\Delta v + \lambda^2 v = 0\mbox{ in }\Omega,\qquad v=1\mbox{ on }\partial\Omega.
\end{equation}
%Exponential substitution
Using the substitution
\begin{equation}\label{eq:exponential}
v(\x)=e^{-\lambda u(\x)}
\end{equation}
in (\ref{eq:screendPoisson}) yields
\begin{equation}\label{eq:heatDist2boundary}
\left(1-|\nabla u|^2\right)+\frac{1}{\lambda}\Delta u=0\mbox{ in }\Omega,\qquad u=0\mbox{ on }\partial\Omega
\end{equation}
and, therefore, since the distance function satisfies the eikonal equation $|\nabla\!\dist(\x,\partial\Omega)|=1$ in $\Omega$, we have
\begin{equation}\label{eq:heatDistApprox}
\dist(\x,\partial\Omega)\approx u(\x)=-\frac{1}{\lambda}\log\left[v(\x)\right],\quad\mbox{as }\lambda\to\infty,
\end{equation}
where $u(\x)$ is the solution to (\ref{eq:heatDist2boundary}).

We can represent $v(\x)$ using the Green's function $G_\lambda(\x)$ for (\ref{eq:screendPoisson}):
$$
v(\x)=\int_{\partial\Omega}\!\frac{\partial G_\lambda(\x-\p)}{\partial\n}\,ds,\quad\x\in\Omega,\quad\p\in\partial\Omega,
$$
where $\n$ is the outer unit normal of $\partial\Omega$. Thus $u(\x)$, the solution to (\ref{eq:heatDist2boundary}), is given by
\begin{equation}\label{eq:GreenDist}
u(\x)\!=\!-\frac{1}{\lambda}\log\!\left[\int_{\partial\Omega}\!\!\!\frac{\partial G_\lambda(\x-\p)}{\partial\n}\,ds\right]\!,
\mbox{ }\p\in\partial\Omega,\mbox{ }\x\in\Omega,
\end{equation}
and one can observe a similarity between (\ref{eq:GreenDist}) and (\ref{eq:logConv}), and also (\ref{eq:underestimate}).

In view of (\ref{eq:dLogConv}), an analogue of (\ref{eq:softMin}) is obtained by differentiating (\ref{eq:GreenDist}) w.r.t. $\lambda$
\begin{equation}\label{eq:dv/v}
-\frac{\partial}{\partial\lambda}\log\left[\int_{\partial\Omega}\!\frac{\partial G_\lambda(\x-\p)}{\partial\n}\,ds\right]=-\frac{\dv(\x)}{v(\x)},
\end{equation}
where $\dv$ denotes the derivative of the solution to (\ref{eq:screendPoisson}) w.r.t. $\lambda$.
In the next section, we will see that, in practice, the right-hand side of (\ref{eq:dv/v}) delivers a reasonably good approximation of the distance function:
\begin{equation}\label{eq:dist1}
\boxed{\dist(\x,\partial\Omega)\approx-\left.\dv(\x)\right/v(\x)}
\end{equation}
Differentiating (\ref{eq:screendPoisson}) by $\lambda$ yields
\begin{equation}\label{eq:screendPoissonD1}
-\Delta\dv+\lambda^2\dv+2\lambda v = 0\quad\mbox{in }\Omega,\quad\dv=0\mbox{ on }\partial\Omega
\end{equation}
and $\dv$ can be found by solving (\ref{eq:screendPoissonD1}).

Now we present a more rigorous way to obtain (\ref{eq:dist1}). Namely let us consider (\ref{eq:heatDist2boundary}) where $s=1/\lambda$ is considered as a small parameter and $u(\x)$, the solution to (\ref{eq:heatDist2boundary}), depends on $s$. One can expect that the first-order Taylor series based estimation of $\left.u(\x)\right|_{s=0}$ delivers a good approximation of the distance
\begin{equation}\label{eq:Taylor1}
\dist(\x,\partial\Omega)=\left.u(\x)\right|_{s=0}\approx u(\x)-s\us(\x),
\end{equation}
where $u(\x)$ and $\us(\x)$ are evaluated at some $s$ close to zero. Simple calculations show that the right-hand sides of (\ref{eq:dist1}) and (\ref{eq:Taylor1}) coincide.

Nothing prevents us from using the second-order Taylor series approximation of $\left.u(\x)\right|_{s=0}$
\begin{equation}\label{eq:Taylor2}
\dist(\x,\partial\Omega)=\left.u(\x)\right|_{s=0}\approx u(\x)-s\us(\x)+\frac{s^2}{2}\uss(\x)
\end{equation}
which yields a possible improvement of (\ref{eq:dist1})
\begin{equation}\label{eq:dist2}
\boxed{\dist(\x,\partial\Omega)\approx-\frac{\dv}{v}
-\frac{\lambda}{2}\left[\frac{\ddv}{v}-\left(\frac{\dv}{v}\right)^2\right]}
\end{equation}
where $\ddv(\x)$ is the solution to the equation obtained from (\ref{eq:screendPoisson}) by differentiating by $\lambda$ twice
$$
-\Delta\ddv+\lambda^2\ddv+2v+4\lambda\dv=0\quad\mbox{in }\Omega,\quad\ddv=0\mbox{ on }\partial\Omega.
$$

If we assume that $u(\x)$ in (\ref{eq:heatDist2boundary}) does not depend on $t$ (this means that $u(\x)$ is the exact distance function), then, in view of (\ref{eq:exponential}), formulas (\ref{eq:heatDistApprox}), (\ref{eq:dist1}), and (\ref{eq:dist2}) are equivalent. In other words, the closer $u(\x)$ is to the exact distance function, the smaller the differences between approximations (\ref{eq:heatDistApprox}), (\ref{eq:Taylor1}), and (\ref{eq:Taylor2}) are.

The original heat method \cite{Crane-Weischedel-Wardetzky_tog13} includes a gradient normalization step which, in our case of estimating $\dist(\x\partial\Omega$, can be described as follows. Given an approximation of the distance function $w(\x)$, we normalize its gradient $g(\x)=\nabla w/|\nabla w|$ and then use the normalized gradient to construct an improved distance function approximation $w_n(\x)$ by solving a simple variational problem
$$
\int_\Omega|\nabla w_n-g|^2d\x\rightarrow\min,\quad w_n(\x)=0\mbox{ on }\partial\Omega.
$$
This least square gradient fitting problem is equivalent to the Poisson problem
$$
\Delta w_n=\div(g)\mbox{ in }\Omega,\quad w_n=0\mbox{ on }\partial\Omega,
$$
which can be efficiently solved numerically. In the next section, we demonstrate that the gradient normalization step improves the approximation accuracy of the right-hand sides of (\ref{eq:dist1}) and (\ref{eq:dist2}).

%%%%%%%%%%%%%%%%%%%%%%%%%%%%%%%%%%%%%%%%%%%%%%%%%%%%%%%%%%%%%%%%%%
% Following the IEEEtran's doc
% https://bigdataieee.org/BigData2020/files/IEEEtran_HOWTO.pdf
% p. 10
\begin{figure*}[htbp]
\centering
\subfloat[]{
\begin{tabular}[b]{c}
\includegraphics[width=.15\textwidth]{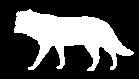}\\
\includegraphics[width=.15\textwidth]{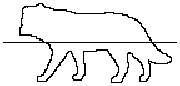}
\end{tabular}
}
\subfloat[]{
\includegraphics[width=.3\textwidth]{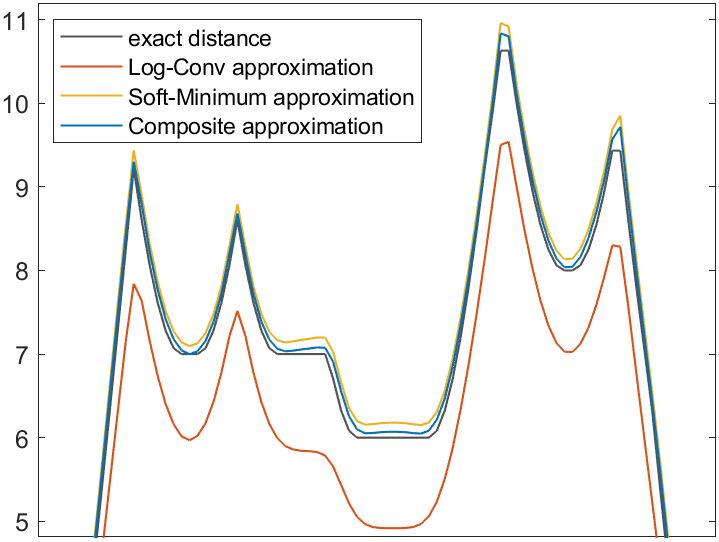}
}
\subfloat[]{
\includegraphics[width=.3\textwidth]{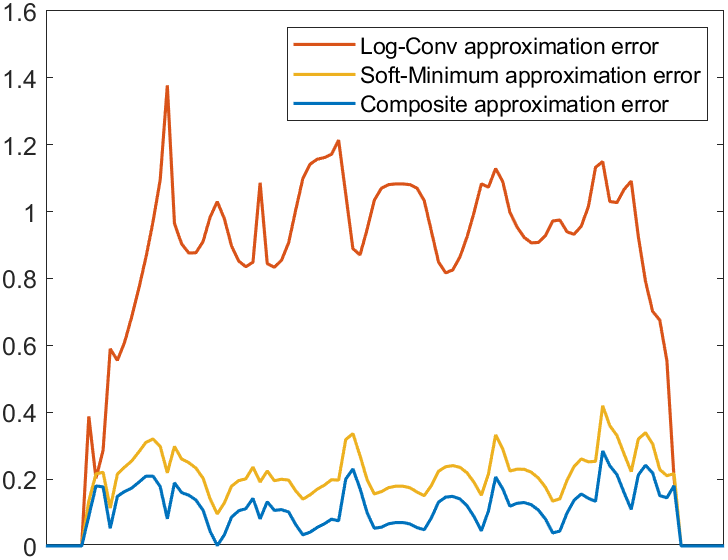}
}
\subfloat[]{
\begin{tabular}[b]{c}
\includegraphics[width=.15\textwidth]{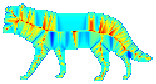}\\
\includegraphics[width=.15\textwidth]{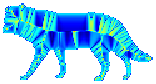}
\end{tabular}
}
\caption{(a) A 2-D shape (binary image) and its 1-D slice which are used to evaluate the proposed combination of the LogConv and SoftMin approximations of the distance to the boundary function. (b) Graphs of the restrictions of the exact distance function, the two convolution-based approximations proposed in \cite{Karam-etal_spl19} and corresponding to the LogConv (\ref{eq:logConv}) and SoftMin (\ref{eq:softMin}) formulas, and the combination of the approximations with weights given by (\ref{eq:weights}) onto the slice. (c) Graphs of the approximation errors.
(d) Approximation errors for the SoftMin (top) and proposed approximation (bottom) are visualized for the whole shape.}
\label{fig:conv2}
\end{figure*}

%%%%%%%%%%%%%%%%%%%%%%%%%%%%%%%%%%%%%%%%%%%%%%%%%%%%%%%%%%%%%%%%%%

\begin{figure*}[t]
\centering
\begin{tabular}{ccc}
\includegraphics[height=4.5cm]{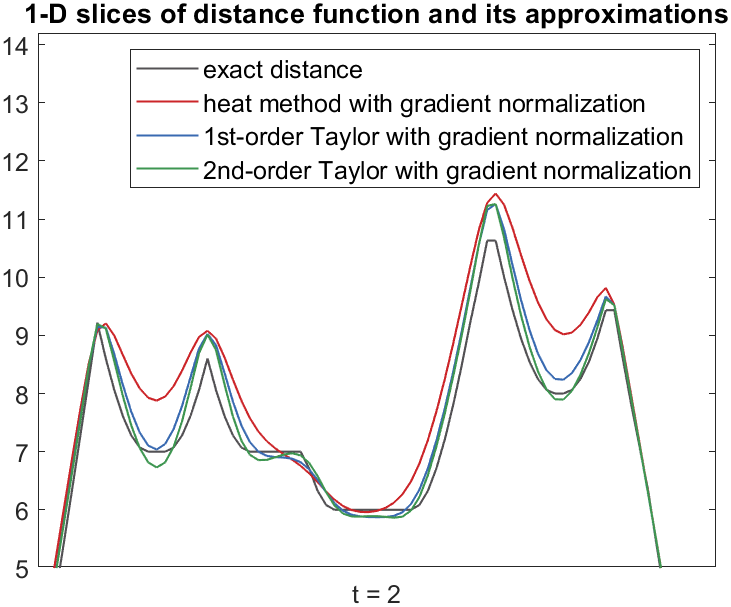}
&
\includegraphics[height=4.5cm]{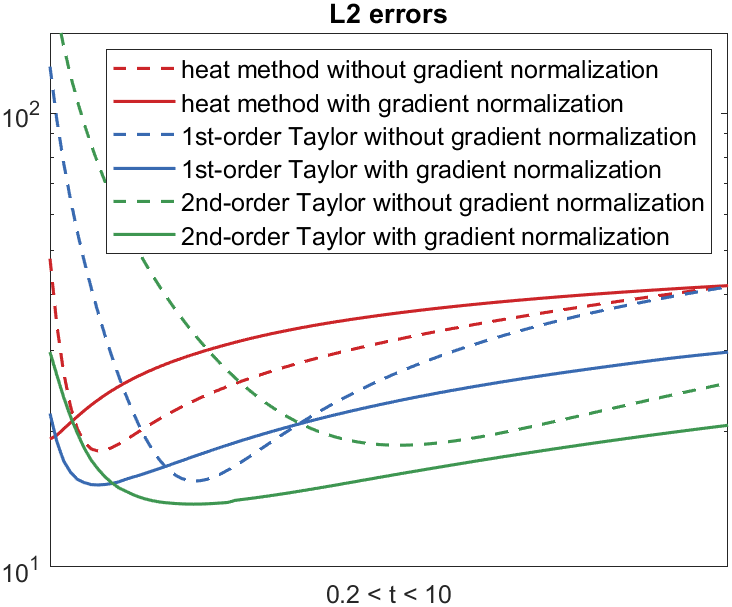}
&
\includegraphics[height=4.5cm]{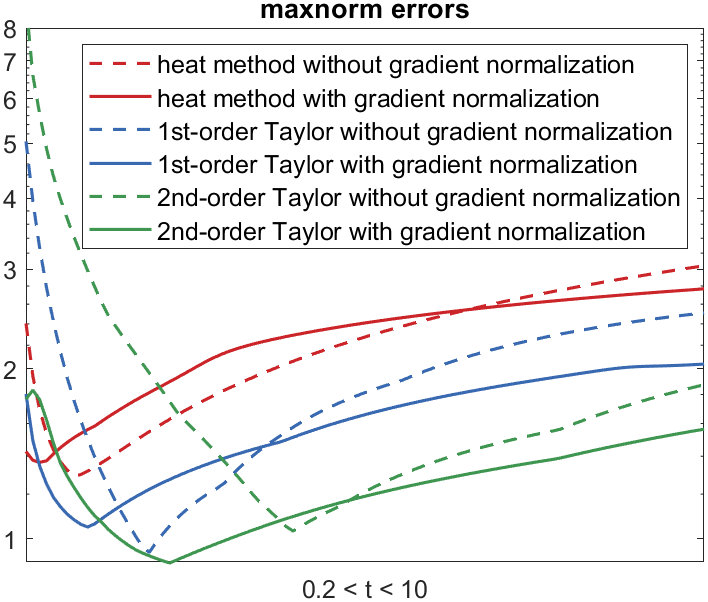}
%\\
%(a) & (b) & (c)
\end{tabular}
\caption{Left: 1-D slices of the distance function approximations (\ref{eq:heatDistApprox}), (\ref{eq:dist1}), and (\ref{eq:dist2}). Middle: dependencies of the $L^2$ distance function approximation errors for (\ref{eq:heatDistApprox}), (\ref{eq:dist1}), and (\ref{eq:dist2}) as functions of parameter $t\equiv1/\lambda^2$. Right: dependencies of the $L^\infty$ (maxnorm) distance function approximation errors for (\ref{eq:heatDistApprox}), (\ref{eq:dist1}), and (\ref{eq:dist2}) as functions of parameter $t$. In all these experiments, the normalized version of (\ref{eq:dist2}) shows the best performance.}
\label{fig:normalization}
\end{figure*}

%%%%%%%%%%%%%%%%%%%%%%%%%%%%%%%%%%%%%%%%%%%%%%%%%%%%%%%%%%%%%%%%%%

\begin{figure*}[h!]
\centering
\begin{tabular}{ccccccc}
\includegraphics[height=4.0cm]{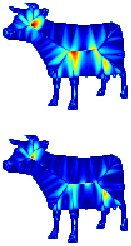}
&
\hspace*{-0.2cm}
\includegraphics[height=4.0cm]{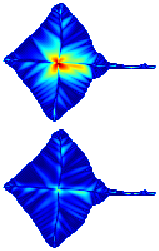}
&
\hspace*{-0.2cm}
\includegraphics[height=4.0cm]{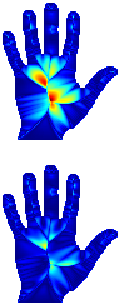}
&
\hspace*{-0.2cm}
\includegraphics[height=4.0cm]{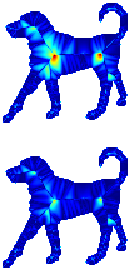}
&
\hspace*{-0.2cm}
\includegraphics[height=4.0cm]{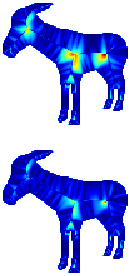}
&
\hspace*{-0.2cm}
\includegraphics[height=4.0cm]{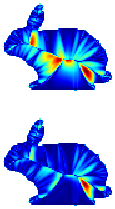}
&
\hspace*{-0.2cm}
\includegraphics[height=4.0cm]{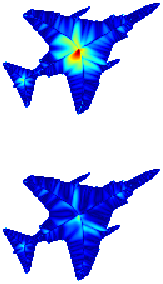}
\\
\includegraphics[height=2.0cm]{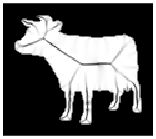}
&
\hspace*{-0.2cm}
\includegraphics[height=2.0cm]{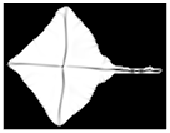}
&
\hspace*{-0.2cm}
\includegraphics[height=2.0cm]{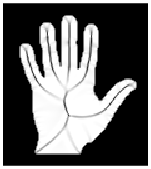}
&
\hspace*{-0.2cm}
\includegraphics[height=2.0cm]{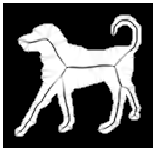}
&
\hspace*{-0.2cm}
\includegraphics[height=2.0cm]{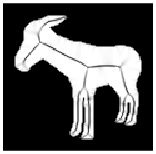}
&
\hspace*{-0.2cm}
\includegraphics[height=2.0cm]{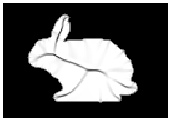}
&
\hspace*{-0.2cm}
\includegraphics[height=2.0cm]{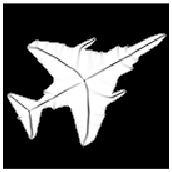}
\end{tabular}
\caption{Distance function error maps for several binary shapes (binary images). Top row: the heat method \cite{Crane-Weischedel-Wardetzky_tog13} with $t=1$ is used. Middle row: second-order Taylor-based extrapolation (\ref{eq:dist2}) with $t\equiv1/\lambda^2=5$ is used. Bottom row: for each shape, the magnitude of the distance function gradient is visualized: high gradient magnitudes correspond to the skeleton of the shape; it looks plausible that the regions near the skeleton branch points contribute most to the distance estimation error.}
\label{fig:errors}
\end{figure*}

%%%%%%%%%%%%%%%%%%%%%%%%%%%%%%%%%%%%%%%%%%%%%%%%%%%%%%%%%%%%%%%%%%
\section{Numerical experiments}
{\em Convolutional approximations.} We start this section by inspecting the idea of blending the right-hand sides of (\ref{eq:logConv}) and (\ref{eq:softMin}) with the weights defined by (\ref{eq:weights}), where constant $K$ is heuristically set equal to $0.1$. The top image of Fig.\,\ref{fig:conv2}(a) presents a 2-D shape (binary image) for which we estimate the distance to the boundary function. For visualization purposes, we consider a one-dimensional slice of the shape, as shown in the bottom image. Fig.\,\ref{fig:conv2}(b) demonstrates graphs of the exact distance function, two its convolution-based approximations proposed in \cite{Karam-etal_spl19} and corresponding to the LogConv (\ref{eq:logConv}) and SoftMin (\ref{eq:softMin}) formulas, and the combination of the approximations with weights given by (\ref{eq:weights}). Fig.\,\ref{fig:conv2}(c) displays the approximation errors. As expected, the combination improves the approximation accuracy. This is also demonstrated by Fig.\,\ref{fig:conv2}(d) where the approximation errors for the SoftMin (top) and proposed approximation (bottom) are visualized for the whole shape (we use a heat color map, where blue corresponds to low errors and red corresponds to high ones).

{\em Differential approximations.}
For the differential distance function approximations (\ref{eq:dist1}) and (\ref{eq:dist2}),
we start by inspecting the approximation properties of their normalized versions. Our main competitor is the heat method \cite{Crane-Weischedel-Wardetzky_tog13,Crane-Weischedel-Wardetzky_cacm17}
and, to match the heat method parameter use, we set $\lambda=1/\sqrt{t}$ in the experiments described below.

In Fig.\,\ref{fig:normalization} we deal with the same shape (binary image) and the same slice of the shape as in Fig.\,\ref{fig:conv2}. The left image of Fig.\,\ref{fig:normalization} presents 1-D slices of the distance function approximations (\ref{eq:heatDistApprox}), (\ref{eq:dist1}), and (\ref{eq:dist2}) with the gradient normalization procedure used for all three approximations. We set $t\equiv1/\lambda^2=2$ in order to achieve better visualization effects. In this example, the difference between (\ref{eq:dist1}), and (\ref{eq:dist2}) is rather small, but both these approximations clearly outperform (\ref{eq:heatDistApprox}).

The middle and right images of Fig.\,\ref{fig:normalization} present the $L^2$ and $L^\infty$ distance function approximation errors for (\ref{eq:heatDistApprox}), (\ref{eq:dist1}), and (\ref{eq:dist2}). The errors are computed for $0.2\leq t\equiv1/\lambda^2\leq10$ for both the normalized and unnormalized versions of the approximations. The normalized version of (\ref{eq:dist2}) demonstrates the best performance while the normalized version of (\ref{eq:dist1}) comes second. One can also observe that the gradient normalization significantly improves the performance of both the schemes (\ref{eq:dist1}), and (\ref{eq:dist2}) proposed in this paper. It is also worth to note that the $L^2$ error for (\ref{eq:dist2}) has a rather flat minimum and, therefore, (\ref{eq:dist2}) demonstrates more or less similar performance for a wide range of values of $t\equiv1/\lambda^2$.

Fig.\,\ref{fig:errors} visualizes distance function error maps computed for the heat method approximation (\ref{eq:heatDistApprox}) and second-order Taylor extrapolation (\ref{eq:dist2}) for a number of binary shapes (images), where hotter colors are used to indicate higher errors and, in particular, red is used to for large errors and blue corresponds to small ones. Normalization is used for both the approximations. We set $t=1$ for (\ref{eq:heatDistApprox}) in order to achieve a high accuracy and avoid possible numerical instability and $t\equiv1/\lambda^2=5$ for (\ref{eq:dist2}). Our choice of the parameter value for (\ref{eq:dist2}) is not optimal but our experiments and the middle image of (\ref{fig:normalization}) suggest that (\ref{eq:dist2}) demonstrates more or less the same performance for a wide range of parameter $t$. As expected, (\ref{eq:dist2}) outperforms (\ref{eq:heatDistApprox}). The images of the bottom row of Fig.\,\ref{fig:errors} suggest that the regions near the skeleton branch points contribute most to the distance estimation error.

%%%%%%%%%%%%%%%%%%%%%%%%%%%%%%%%%%%%%%%%%%%%%%%%%%%%%%%%%%%%%%%%%%
\section{Discussion, conclusion, and future work}
Our Matlab implementations of (\ref{eq:dist1}) and (\ref{eq:dist2}) and corresponding distance function approximations will be made available upon acceptance of this letter. The implementations are simple and can be easily obtained by introducing proper modifications into the Matlab scripts presented in \cite{Belyaev-Fayolle_na20}.

The main limitation of this paper is that we deal with the planar distance estimation problem. It has certain practical applications including its use within the so-called characteristic function method \cite{Babuska-etal_an03}. However the more challenging curvilinear and multidimensional distance function estimation problems are not considered in the paper. Another limitation is a lack of mathematical justification of the proposed Taylor-based extrapolation schemes (\ref{eq:dist1}) and (\ref{eq:dist2}). They look plausible but in contrast to the heat method \cite{Crane-Weischedel-Wardetzky_tog13} which is based on mathematical results of S.\,R.\,S.\,Varadhan \cite{Varadhan_cpam67} are not supported by truly rigorous mathematical results.

It is interesting to observe a certain similarity between the convolutional and differential distance function estimation schemes considered in this paper and a PDE-based surface reconstruction approach developed in \cite{Kazhdan_sgp05,Kazhdan_etal_sgp06,Kazhdan-Hoppe_tog13}, where the evolution from convolutional to PDE-based methods is, in some sense, similar to how we move from the convolutional schemes to the differential ones,

It would be interesting to check if the convolution-based distance function approximation schemes we introduced in this paper can be extended to graphs by using graph Laplacian convolutions. We hope that, similar to the heat method \cite{Crane-Weischedel-Wardetzky_tog13}, PDE-based distance function approximations (\ref{eq:Taylor1}) and (\ref{eq:Taylor2}) can be used for computing geodesic distances on surfaces and, in particular, we would like to check if our approximations can compete with \cite{Feng-Crane_tog24} in terms of robustness. We are also interested in extending (\ref{eq:Taylor1}) and (\ref{eq:Taylor2}) to Finsler metrics \cite{Yang-etal_accv19,Weber-etal_cvpr24}.
Adapting our distance function approximations for path planning applications \cite{Muchacho-Pokorny_arXiv24} constitute another direction for future research.
Finally, it would be interesting to see if Pad\'{e} approximations are capable of outperforming our Taylor-based extrapolation approach.
%%%%%%%%%%%%%%%%%%%%%%%%%%%%%%%%%%%%%%%%%%%%%%%%%%%%%%%%%%%%%%%%%
\newpage\phantom{m}\newpage\phantom{m}\newpage
\bibliographystyle{IEEEtran}
\bibliography{distBib}

\begin{thebibliography}{10}
\providecommand{\url}[1]{#1}
\csname url@samestyle\endcsname
\providecommand{\newblock}{\relax}
\providecommand{\bibinfo}[2]{#2}
\providecommand{\BIBentrySTDinterwordspacing}{\spaceskip=0pt\relax}
\providecommand{\BIBentryALTinterwordstretchfactor}{4}
\providecommand{\BIBentryALTinterwordspacing}{\spaceskip=\fontdimen2\font plus
\BIBentryALTinterwordstretchfactor\fontdimen3\font minus
  \fontdimen4\font\relax}
\providecommand{\BIBforeignlanguage}[2]{{%
\expandafter\ifx\csname l@#1\endcsname\relax
\typeout{** WARNING: IEEEtran.bst: No hyphenation pattern has been}%
\typeout{** loaded for the language `#1'. Using the pattern for}%
\typeout{** the default language instead.}%
\else
\language=\csname l@#1\endcsname
\fi
#2}}
\providecommand{\BIBdecl}{\relax}
\BIBdecl

\bibitem{Park-etal_jcp24}
Y.~Park, C.~h. Song, J.~Hahn, and M.~Kang, ``{ReSDF:} {R}edistancing implicit
  surfaces using neural networks,'' \emph{Journal of Computational Physics},
  vol. 502, 2024.

\bibitem{King-etal_tog24}
N.~King, H.~Su, M.~Aanjaneya, S.~Ruuth, and C.~Batty, ``A closest point method
  for {PDEs} on manifolds with interior boundary conditions for geometry
  processing,'' \emph{ACM Transactions on Graphics}, 2024.

\bibitem{MattosDaSilva-etal_tog24}
L.~Mattos Da~Silva, O.~Stein, and J.~Solomon, ``A framework for solving
  parabolic partial differential equations on discrete domains,'' \emph{ACM
  Transactions on Graphics}, 2024.

\bibitem{Huguet-etal_NeurIPS23}
G.~Huguet, A.~Tong, E.~De~Brouwer, Y.~Zhang, G.~Wolf, I.~Adelstein, and
  S.~Krishnaswamy, ``A heat diffusion perspective on geodesic preserving
  dimensionality reduction,'' in \emph{37th Conference on Neural Information
  Processing Systems (NeurIPS 2023)}, 2023.

\bibitem{Edelstein-etal_sig23}
M.~Edelstein, N.~Guillen, J.~Solomon, and M.~Ben-Chen, ``A convex optimization
  framework for regularized geodesic distances,'' in \emph{ACM SIGGRAPH 2023
  Conference Proceedings}, 2023, pp. 2:1--2:11.

\bibitem{Feng-Crane_tog24}
N.~Feng and K.~Crane, ``A heat method for generalized signed distance,''
  \emph{ACM Transactions on Graphics}, vol.~43, no.~4, pp. 92:1--92:16, 2024.

\bibitem{Yang-etal_NeurIPS23}
H.~Yang, Y.~Sun, G.~Sundaramoorthi, and A.~Yezzi, ``Stabilizing the
  optimization of neural signed distance functions and finer shape
  representation,'' in \emph{37th Conference on Neural Information Processing
  Systems (NeurIPS 2023)}, 2023.

\bibitem{Meister_arXiv23}
M.~Meister, ``A fast algorithm for all-pairs-shortest-paths suitable for neural
  networks,'' \emph{arXiv preprint arXiv:2308.07403}, 2023.

\bibitem{Karam-etal_spl19}
C.~Karam, K.~Sugimoto, and K.~Hirakawa, ``Fast convolutional distance
  transform,'' \emph{IEEE Signal Processing Letters}, vol.~26, no.~6, pp.
  853--857, 2019.

\bibitem{Tibshirani-etal_arXiv24}
R.~J. Tibshirani, S.~W. Fung, H.~Heaton, and S.~Osher, ``{L}aplace meets
  {M}oreau: {S}mooth approximation to infimal convolutions using {L}aplace's
  method,'' \emph{arXiv preprint arXiv:2406.02003}, 2024.

\bibitem{Varadhan_cpam67}
S.~R.~S. Varadhan, ``On the behavior of the fundamental solution of the heat
  equation with variable coefficients,'' \emph{Comm. Pure Appl. Math.},
  vol.~20, pp. 431--455, 1967.

\bibitem{Crane-Weischedel-Wardetzky_tog13}
K.~Crane, C.~Weischedel, and M.~Wardetzky, ``Geodesics in heat: A new approach
  to computing distance based on heat flow,'' \emph{ACM Transactions on
  Graphics}, vol.~32, pp. 152:1--152:11, 2013.

\bibitem{Crane-Weischedel-Wardetzky_cacm17}
------, ``The heat method for distance computation,'' \emph{Communications of
  the ACM}, vol.~60, no. 112, pp. 90--99, 2017.

\bibitem{Belyaev-Fayolle_na20}
A.~Belyaev and P.-A. Fayolle, ``An {ADMM}-based scheme for distance function
  approximation,'' \emph{Numerical Algorithms}, vol.~84, pp. 983--996, 2020.

\bibitem{Babuska-etal_an03}
I.~Babu{\v{s}}ka, U.~Banerjee, and J.~E. Osborn, ``Survey of meshless and
  generalized finite element methods: {A} unified approach,'' \emph{Acta
  Numerica}, vol.~12, pp. 1--125, 2003.

\bibitem{Kazhdan_sgp05}
M.~Kazhdan, ``Reconstruction of solid models from oriented point sets,'' in
  \emph{Proceedings of the third Eurographics symposium on Geometry
  processing}, 2005, pp. 73--82.

\bibitem{Kazhdan_etal_sgp06}
M.~Kazhdan, M.~Bolitho, and H.~Hoppe, ``Poisson surface reconstruction,'' in
  \emph{Proceedings of the fourth Eurographics symposium on Geometry
  processing}, vol.~7, no.~4, 2006.

\bibitem{Kazhdan-Hoppe_tog13}
M.~Kazhdan and H.~Hoppe, ``Screened {P}oisson surface reconstruction,''
  \emph{ACM Transactions on Graphics}, vol.~32, no.~3, pp. 1--13, 2013.

\bibitem{Yang-etal_accv19}
F.~Yang, L.~Chai, D.~Chen, and L.~Cohen, ``Geodesic via asymmetric heat
  diffusion based on {F}insler metric,'' in \emph{Computer Vision--ACCV 2018:
  14th Asian Conference on Computer Vision, Revised Selected Papers, Part V
  14}, 2019, pp. 371--386.

\bibitem{Weber-etal_cvpr24}
S.~Weber, T.~Dag{\`e}s, M.~Gao, and D.~Cremers,
  ``{F}insler-{L}aplace-{B}eltrami operators with application to shape
  analysis,'' in \emph{Proceedings of the IEEE/CVF Conference on Computer
  Vision and Pattern Recognition}, 2024, pp. 3131--3140.

\bibitem{Muchacho-Pokorny_arXiv24}
R.~I.~C. Muchacho and F.~T. Pokorny, ``Walk on spheres for {PDE}-based path
  planning,'' \emph{arXiv preprint arXiv:2406.01713}, 2024.

\end{thebibliography}
%%%%%%%%%%%%%%%%%%%%%%%%%%%%%%%%%%%%%%%%%%%%%%%%%%%%%%%%%%%%%%%%%%
\end{document}